\documentclass[12pt,a4paper,notitlepage]{article}
\usepackage{amsmath}
\usepackage{graphicx}
\makeatletter
\usepackage{hhline}
\usepackage{amsthm}
\usepackage{amscd}
\usepackage{amssymb}
\usepackage{latexsym}
\usepackage{tabls}

\newtheorem{theorem}{Theorem}[section]
\newtheorem{proposition}[theorem]{Proposition}
\newtheorem{lemma}[theorem]{Lemma}

\newtheorem{example}[theorem]{Example}
\newtheorem{definition}[theorem]{Definition}
\newtheorem{remark}[theorem]{Remark}
\newtheorem*{acknowledge}{Acknowledgments}

\def\cosec{\operatorname{cosec}}

\def\Ker{\operatorname{Ker}}
\def\loc{\operatorname{loc}}
\def\mon{\operatorname{mon}}
\def\Sign{\operatorname{Sign}}

\def\Int{\operatorname{Int}}
\def\Diff{\operatorname{Diff}}

\def\Deck{\operatorname{Deck}}

\def\mod{\operatorname{mod}}

\begin{document}
\title{A local signature for fibrations with a finite group action}
\author{Masatoshi Sato}
\date{\empty}
\maketitle
\begin{abstract}
Let $p$ be a finite regular covering on a 2-sphere with at least three branch points. In this paper, we construct a local signature for the class of fibrations whose general fibers are isomorphic to the covering $p$.
\end{abstract}

\section{Introduction}
Let $E$ and $B$ be closed oriented smooth manifolds of dimension 4 and 2, respectively. Let $\Sigma_g$ denote a closed oriented surface of genus $g\ge1$. Assume that a smooth surjective map $f:E\to B$ has finitely many critical values $\{b_l\}_{l=1}^n$, and the fiber $f^{-1}(b)$ on $b\in B-\{b_l\}_{l=1}^n$ is connected. Then, its restriction $E-\{f^{-1}(b_l)\}_{l=1}^n \to B-\{b_l\}_{l=1}^n$ is an oriented fiber bundle whose fiber is diffeomorphic to $\Sigma_g$. We call the triple $(f,E,B)$ satisfying these conditions a topological fibration of genus $g$. Inverse images of a regular value and a singular value under $f$ are called a regular fiber and a singular fiber, respectively.

For a fibration $f:E\to B$, denote by $\Delta_l\subset B$ a closed neighborhood of the critical value $b_l$. Denote by $E_l$ the inverse image $f^{-1}(\Delta_l)$, and the restriction $f|_{E_l}$ by $f_l$. On some classes of fibrations, the signature of the fibration $f:E\to B$ is described as the sum
\[
\Sign E=\sum_{l=1}^n\sigma_{\loc}(f_l,E_l,\Delta_l),
\]
of local invariants $\sigma_{\loc}(f_l,E_l,\Delta_l)\in\mathbf{R}$ each of which depends on a neighborhood of a singular fiber. To be precise, these local invariants are defined as a function $\sigma_{\loc}$ on the set of singular fiber germs which arise in the class of fibrations. We call this function the local signature.

One of the motivations for the study of local signatures is that it is closely related to a cobounding function of the Meyer cocycle, an important 2-cocycle of the mapping class group of the surface $\Sigma_g$. This function is related to several invariants including the eta invariant of the signature operator (see Atiyah \cite{atiyah1987ldeta}, Iida \cite{iida2004ale}) and the Casson invariant of homology 3-spheres (see Morita \cite{morita1989csi}). In algebraic geometry, the local signature is related to the slope equality problem. This is studied in order to describe the geography of algebraic surfaces of general type (see Ashikaga-Konno \cite{ashikaga2000gal}). 

For many classes of fibrations, local signatures are constructed and calculated in various fields including topology, algebraic geometry, and complex analysis. For the fibrations of genus 1 and 2, Matsumoto \cite{matsumoto1983mft}\cite{matsumoto1996lfg} constructed local signatures using the cobounding function of the Meyer cocycle. He also calculated the values for Lefschetz singular fiber germs. Ueno also constructed and calculated the local signature for the genus 2 fibrations using the even theta constant. For the fibrations of genus more than 2 whose monodromies are in the whole mapping class group, there does not exist a local signature. But for some restricted class, such as hyperelliptic fibrations, Endo \cite{endo2000mss} constructed it (see also Morifuji \cite{morifuji2003msf}). Arakawa and Ashikaga also constructed it in the setting of algebraic geometry, and Terasoma proved that these two local signatures coincides. Local signatures for many kinds of restricted classes of fibrations are listed in Ashikaga-Endo \cite{ashikaga2006vad} and Ashikaga-Konno \cite{ashikaga2000gal}. Recently, Iida \cite{iida2004ale}, Kuno \cite{kuno2007mcg}\cite{kuno2009mfp} and Yoshikawa \cite{yoshikawa2000lsg} constructed it for some restricted classes.

The purpose of this paper is to construct the local signatures for the classes of topological fibrations which have a fiber-preserving finite group action. We also assume that the quotient space of their general fiber is a sphere with at least 3 branch points. We should mention that Furuta \cite{furuta1999lsbse} constructed a local signature for broader classes of fibrations than ours in the manner of differential geometry. Furthermore, Nakata \cite{nakata2005lsle} calculated it on Lefschetz fiber germs of the hyperelliptic fibration. But the local signature in this paper is easier to compute than that one. In general, a local signature of this class is not unique. I do not know whether the local signature in this paper coincides with Furuta's one.

Let $G$ be a finite group, and $\Sigma$ a closed surface. We call a finite regular covering $p:\Sigma_g\to \Sigma$ a $G$-covering if its deck transformation group is isomorphic $G$. For a $G$-covering $p$, Birman-Hilden \cite{birman1973ihr} defined a group called the symmetric mapping class group. The monodromy group of fibrations of the $G$-covering $p$ is considered as the subgroup of this group. The local signature induces a cobounding function of the pullback of the Meyer cocycle in the symmetric mapping class group. Using this cobounding function, we will construct the local signature which can be applied to a broader class of fibrations than our class.  

This paper is organized as follows. In Section 2, we define a class of fibrations for which we will construct the local signature later. We also define a broader class of fibrations and review the local signature for it constructed by Furuta. For a fibration $f:E\to B$ of the $G$-covering $p$, let $E^h\subset E$ denote the fixed point set in $E$ for $h\in G$. Since general fibers of the fibration are isomorphic to the $G$-covering $p$, the subspace $E^h\cup(E-\{f^{-1}(b_l)\}_{l=1}^n)$ of $E^h$ can be considered as a (not necessarily connected) covering space of $B-\{b_l\}_{l=1}^n$. We call a component $S$ of $E^h$ horizontal if $S\cup(E-\{f^{-1}(b_l)\}_{l=1}^n)$ is a covering space of $B-\{b_l\}_{l=1}^n$. The fixed point set $E^h$ consists of these horizontal components and vertical components which are included in singular fibers. For $\psi\in[0,\pi]$, denote by $I(h,\psi)$ the union of the horizontal components whose normal bundles rotate $\pm\psi$ under the action of $h\in G$. Let $J(f_l,h)$ denote the union of vertical components included in a fiber germ $f^{-1}(b_l)$. In Section 3, we prove the main theorem (Theorem 1.1) as below assuming that the sum of normal euler numbers $\chi(N(S))$ of the connected component $S$ of $I(h,\psi)$ for $h\in G$ localizes as
\[
\sum_{S\subset I(h,\psi)}\chi(N(S))=\sum_{l=1}^n \chi_{\loc}^{h,\psi}([f_l,E_l,\Delta_l]).
\]
Here $\chi_{\loc}^{h,\psi}([f_l,E_l,\Delta_l])$ is a rational number which depends on the fiber germ $[f_l,E_l,\Delta_l]$ in the fibration $f:E\to B$.

Let $f:E\to \Delta$ be a fiber germ of the $G$-covering $p$ on a closed 2-disk $\Delta$. Let $\{P_j\}$ and $\{F_j\}$ be the discrete fixed points and 2-dimensional components in $J(f,h)$. Denote the rotating angle of the normal bundle of $F_j$ by $\pm\psi_j$ for $\psi_j\in[0,\pi]$. Identifying the neighborhood of $P_j$ and $\mathbf{C}^2$, choose $\varphi_j$, $\varphi_j'\in[0,2\pi]$ so that the action of $h$ is written as $(z,w)\mapsto (e^{\sqrt{-1}\varphi_j}z,e^{\sqrt{-1}\varphi_j'}w)$ in a suitable coordinate. With the local normal euler number $\chi_{\loc}^{h,\psi}$ and these  connected components of the fixed point set, The local signature $\sigma_{\loc}$ is described as follows. 
\begin{theorem}\label{thm:main}
Let $p:\Sigma_g\to S^2$ be a $G$-covering with at least three branch points. The signatures of fibrations of the $G$-covering $p$ localizes. Our local signature is written as
\begin{align*}
\sigma_{\loc}([f,E,\Delta])=|G|\Sign&(E/G)\\
+\sum_{h\ne1\in G}&\bigg(-\sum_{\psi\in [0,\pi]}\chi_{\loc}^{h,\psi}([f,E,\Delta])\cosec^2\left(\frac{\psi}{2}\right)\\
&+\sum_{P_j\subset J(f,h)}\cot\left(\frac{\varphi_j}{2}\right)\cot\left(\frac{\varphi_j\smash{'}}{2}\right)-\sum_{F_j\subset J(f,h)}\chi(N(F_j))\cosec^2\left(\frac{\psi_j}{2}\right)\bigg).
\end{align*}
\end{theorem}
The key tools to prove this theorem are the localization of the normal euler number and the $G$-signature theorem (Atiyah-Singer \cite{atiyah1968ieo}). In Section 4, we construct a local euler number using the multi-section on the normal bundle of $I(h,\psi)$ made by Furuta, and complete the proof of the main theorem (Theorem \ref{thm:main}). 

In the rest of paper, we will consider a cobounding function of the pullback of the Meyer cocycle in the symmetric mapping class group of the $G$-covering $p$. In Section 5, we construct a local signature for broader class of fibrations of the $G$-covering $p$ when $p$ satisfies some condition. To do this, we describe the local signature using the cobounding function of the pullback of the Meyer cocycle in the symmetric mapping class group. In Section 6, we give a standard generating system of the symmetric mapping class group of a $G$-covering when $G$ is abelian. Let $d\ge2$ and $m\ge3$ be integers such that $m$ is divided by $d$, and $A$ a finite set $\{\alpha_i\}_{i=1}^m$ in $S^2$. For each $i=1,2,\cdots,m$, choose a loop $\gamma_{\alpha_i}$ which rotates around a point $\alpha_i$ counterclockwise once. Define a surjective homomorphism $k:H_{1}(S^{2}-A)\to\mathbf{Z}_{d}$ by mapping each homology class $[\gamma_{\alpha_i}]$ to $1 \mod d$. Let $p_1:\Sigma_g\to S^2$ be the $\mathbf{Z}_d$-covering on $S^2$ which has the branch set $A$ in $S^{2}$ and the monodromy homomorphism $k$. In Section 7, we calculate the local signature of a fiber germ $f:E\to\Delta$ of the $\mathbf{Z}_d$-covering $p_1$. Its monodromy is the inverse of an element $\hat{\sigma}_{ij}$ in the standard generating system of the symmetric mapping class group for the $\mathbf{Z}_d$-covering $p_1$. We also calculate the value of the cobounding function $\phi(\hat{\sigma}_{ij})$.

\begin{proposition}\label{proposition:localsignature}
Let $f:E\to\Delta_1$ be the representative of a fiber germ in $S_g^{p_1}$ constructed in the proof of Lemma \ref{cover p1}. Then we have
\[
\sigma_{\loc}([f,E,\Delta_1])=-\dfrac{(d-1)(d+1)m}{3d(m-1)},
\]
and
\[
\phi(\hat{\sigma}_{ij})=\dfrac{(d-1)(d+1)m}{3d(m-1)}
\]
where $\hat{\sigma}_{ij}$ is the generator of $\mathcal{M}_g(p_1)$ defined in Section \ref{section:generator}.
\end{proposition}

\section{Fibrations of the $G$-covering $p$}\label{section:fibrations}
Denote by $\Sigma_g$ a closed oriented surface of genus $g\ge1$. Let $G$ be a finite group, and $p:\Sigma_g\to S^2$ a $G$-covering, that is, a finite regular covering whose deck transformation group $\Deck(p)$ is isomorphic to $G$. In the sequel, we fix an isomorphism between the deck transformation group $\Deck(p)$ and $G$. In this section, we define two kinds of fibrations of the $G$-covering $p$ in Definition \ref{broad fibration} and Definition \ref{narrow fibration}. Later in Section \ref{euler-localize}, we will construct a local signature for the fibrations in Definition \ref{narrow fibration}. 

Let $C(p)$ denote the centralizer of $\Deck(p)$ in the orientation-preserving diffeomorphism group $\Diff_{+}\Sigma_{g}$ of the surface $\Sigma_g$.

\begin{definition}\label{broad fibration}
Let $E$ and $B$ be compact (not necessarily closed) manifolds of dimension 4 and 2, and $f:E\to B$ a smooth surjective map. A triple $(f,E,B)$ is called a fibration of the $G$-covering $p$ in a broad sense if it satisfies
\begin{enumerate}
\item $\partial E=f^{-1}(\partial B)$,
\item $f:E\to B$ has finitely many critical values $\{b_l\}_{l=1}^n$ in $\Int B$, and the restriction $E-\{f^{-1}(b_l)\}_{l=1}^n \to B-\{b_l\}_{l=1}^n$ is a smooth oriented $\Sigma_g$-bundle,
\item The structure group of the $\Sigma_{g}$-bundle $E-\{f^{-1}(b_{l})\}_{l=1}^{n}\to B-\{b_{l}\}_{l=1}^{n}$ is included in $C(p)$.
\end{enumerate}
\end{definition}

The fibration of the $G$-covering is a generalization of the hyperelliptic fibration, a fibration of a $\mathbf{Z}_{2}$-covering $\Sigma_g\to S^2$ on a sphere. The natural action of $G$ on $\Sigma_g$ gives rise to a smooth fiberwise $G$-action on $E-\{f^{-1}(b_{l})\}_{l=1}^n$, since the structure group of the fiber bundle $E-\{f^{-1}(b_{l})\}_{l=1}^n$ is included in $C(p)$. Note that for each regular value $b\in B$, the covering $f^{-1}(b)\to f^{-1}(b)/G$ is isomorphic to the $G$-covering $p$.

On a class of fibrations of a covering which is not necessarily regular, Furuta has already constructed a local signature. For a more detailed setting, see Furuta \cite{furuta1999lsbse}. He constructed a canonical multi-section on the relative tangent bundle of the fiber bundle $E-\{f^{-1}(b_l)\}_{l=1}^n\to B-\{b_l\}_{l=1}^n$, using the fact that any fiber has at least 3 branch points. He made a connection on the tangent bundle $TE$ by using the multi-section, and showed that its Pontrjagin form vanishes outside neighborhoods of the singular fibers $\{f^{-1}(b_l)\}_{l=1}^n$. Thus, the signature localizes. 

We also use this multi-section to construct our local signature. But, in this paper, we consider a narrower class of fibrations in order to make a local signature without using the connection which is easy to compute. A fibration of the $G$-covering $p$ (in a narrow sense) is defined as follows.
\begin{definition}\label{narrow fibration}
A triple $(f, E, B)$ is called a fibration of the $G$-covering $p$ (in a narrow sense), if it satisfies
\begin{enumerate}
\item The map $f:E\to B$ is a fibration of the $G$-covering $p$ in the broad sense, 
\item The natural $G$-action on $E-\{f^{-1}(b_{l})\}_{l=1}^n$ extends to a smooth action on $E$.
\end{enumerate}
\end{definition}

In this paper, we simply call it a fibration of the $G$-covering $p$. Our local signature and local Euler number is defined as functions on the set of fiber germs which arise in these fibrations. The set of fiber germs is defined as follows. Denote by $\Delta$ a closed 2-disk. Consider fibrations ($f$,$E$,$\Delta$) of the $G$-covering $p$ with unique critical values $b\in\Delta$. Let ($f_1$,$E_1$,$\Delta_1$) and ($f_2$,$E_2$,$\Delta_2$) be such fibrations which have unique critical values $b_1$ and $b_2$, respectively. We call these fibrations are equivalent if and only if there exist closed 2-disks $\Delta_1'\subset\Delta_1$ and $\Delta_2'\subset\Delta_2$ including the critical values, an orientation-preserving diffeomorphism $\varphi:(\Delta_1',b_1)\to(\Delta_2',b_2)$, and a $G$-equivariant orientation-preserving diffeomorphism $\tilde{\varphi}:f_1^{-1}(\Delta_1)\to f_2^{-1}(\Delta_2)$ such that
\[
\varphi f_1=f_2\tilde{\varphi}.
\]
We call this equivalent class a fiber germ of fibrations of the $G$-covering $p$, and denote the set of equivalent classes by $S_g^p$.

\section{A local signature on the class of fibrations of the $G$-covering $p$ in the narrow sense}\label{section:local signature}
Let $G$ be a finite group, and $p:\Sigma_g\to S^2$ a $G$-covering as in Section \ref{section:fibrations}. Let $(f,E,B)$ be a fibration of the $G$-covering $p$. We assume $E$ and $B$ are without boundary. For $h\in G$ and $\psi\in [0,\pi]$, denote by $I(h,\psi)$ the union of all horizontal components of the fixed point set $E^h$ whose normal bundles are rotated $\pm\psi$ by the action of $h\in G$. Let $\{[f_{l}, E_l, \Delta_l]\}_{l=1}^{n}$ be the fiber germs in $f:E\to B$, and denote by $\chi(N(S))$ the normal euler number of a connected component $S$ in $I(h,\psi)$. If there exists a function
\[
\chi_{\loc}^{h,\psi}:S_{g}^{p}\to\mathbf{Q}
\]
on the set of fiber germs such that
\[
\sum_{S\subset I(h,\psi)}\chi(N(S))=\sum_{l=1}^{n}\chi_{\loc}^{h,\psi}([f_{l},E_l,\Delta_l]),
\]
we say that the normal euler number of the horizontal components $I(h,\psi)$ localizes. In this section, we prove the main theorem (Theorem \ref{thm:main}) assuming that this number localizes. We will construct a local euler number later in Definition \ref{local euler number}.

\subsection{The fixed point set of the $G$-action}\label{fixed-point-set}
To construct the local signature, we will apply the $G$-signature theorem on the total space $E$ of the fibration $f:E\to B$ of the $G$-covering $p$. Hence the fixed point set of $G$-action on $E$ plays an important role. We investigate this set in this section.

For $h\in G$, the fixed point set $E^h$ is a pairwise disjoint collection of closed submanifolds (See, for example, Conner \cite{conner1979dpm} p.72). Since the group $G$ preserves the orientation of $E$, it consists of closed $2$-manifolds $\{S_{i}\}$ and 0-manifolds $\{P_{j}\}$.  In Introduction, we defined two kinds of components of $E^h$. A component of $E^h$ is called vertical if they are contained in a singular fiber, and it is called horizontal otherwise. Let $N(b)$ denote a neighborhood of a regular value $b\in B$ of $f$. If we endow the natural $G$-action on the first factor of $\Sigma_g\times D^2$, a local trivialization $f^{-1}(N(b))\cong \Sigma_g\times D^2$ preserving the $C(p)$ structure is $G$-equivariant. Therefore, any horizontal component $S$ is 2-dimensional, and $S\cap f^{-1}(B-\{b_{l}\}_{l=1}^n)$ is a (not necessarily connected) covering space of $B-\{b_{l}\}_{l=1}^n$. We give an example of horizontal components.
\begin{example}[An elliptic surface]
Let $[x:y:z]$ and $[b_{1}:b_{2}]$ be a homogeneous coordinate of $\mathbf{CP}^{2}$ and $\mathbf{CP}^{1}$. Consider the singular surface
\[
E=\{b_{2}y^{2}z=x(x-z)(b_{2}x-b_{1}z)\,|\,([x:y:z],[b_{1}:b_{2}])\in\mathbf{CP}^{2}\times\mathbf{CP}^{1}\}.
\]
If we blow up $E$ at 
\[
([x:y:z],[b_1:b_2])=([0:0:1],[0,1]),([1:0:1],[1,1]),([1:0:0],[1,0]),
\]
we obtain a smooth elliptic surface $E'$ and a smooth map $E'\to \mathbf{CP}^1$ defined by the natural projection. Endow the action of $\mathbf{Z}_{2}$ on $E$ by
\[
t:([x:y:z],[b_{1}:b_{2}])\mapsto([x:-y:z],[b_{1}:b_{2}]).
\]
This group action extends naturally onto $E'$, and its fixed point set is the disjoint sum of the proper transforms of $\{x=y=0\}$, $\{x=z,y=0\}$, $\{b_{2}x=b_{1}z,y=0\}$, and $\{x=z=0\}$. These components are horizontal.
\end{example}

\subsection{Localizations of signature and the normal euler number}\label{euler-localize}
We deduce the following lemma from the $G$-Signature theorem.
\begin{lemma}\label{localize}
Let $f:E\to B$ a fibration of the $G$-covering $p$. If the normal euler number of the horizontal components $I(h,\psi)$ localizes for any $h\in G$ and $\psi\in [0,\pi]$, the signature of $E$ also localizes.
\end{lemma}
To prove this lemma, we prepare some facts about the $G$-signature. Let $G$ be a finite group which acts on a closed oriented 4-manifold $X$ preserving the orientation. Then, the $G$-signature $\Sign(h,X)$ for $h\in G$ is a rational number satisfying 
\begin{equation}
\Sign X=-\sum_{h\ne1\in G}\Sign(h,X)+|G|\Sign(X/G).\label{aver}
\end{equation}
For the details, see Atiyah-Singer \cite{atiyah1968ieo} (See also Gordon \cite{gordon1986gst}).

The fixed point set of $h\in G$ is the disjoint sum of closed 2-manifolds $\{S_{i}\}$ and 0-manifolds $\{P_{j}\}$. Denote the rotation angle of the normal bundle of $S_{i}$ by $\pm\psi_{i}$, where $\psi\in [0,\pi]$. Identifying neighborhoods of $P_j$ and the origin in $\mathbf{C}^2$, choose $\varphi_j$, $\varphi_j'\in[0,2\pi]$ so that the action of $h$ is written as $(z,w)\mapsto (e^{\sqrt{-1}\varphi_j}z,e^{\sqrt{-1}\varphi_j'}w)$ in a suitable coordinate. Then, the $G$-signature $\Sign(h,X)$ is written in terms of these rotation angles and the normal euler number of $S_{i}$ as follows.
\begin{theorem}[Atiyah-Singer \cite{atiyah1968ieo}]
\begin{equation}
\Sign(h,X)=\sum_{i}\chi(N(S_{i}))\cosec^{2}(\frac{\psi_{i}}{2})-\sum_{j}\cot(\frac{\varphi_{j}}{2})\cot(\frac{\varphi'_{j}}{2}).\label{G-sign}
\end{equation}
\end{theorem}

\begin{proof}[proof of lemma \ref{localize}]
Choose representatives $\{f_l,E_l,\Delta_l\}_{l=1}^n$ of all fiber germs in the fibration $f:E\to B$. Note that the complement of the fiber germs $\bigcup_{l=1}^n E_l/G$ in the quotient space $E/G$ is a $S^2$-bundle. Since the signature of a $S^2$-bundle on a compact 2-manifold vanishes, we have
\begin{align}
\Sign(E/G) & =\Sign((E-\amalg_{l=1}^n \Int E_{l})/G)+\sum_{l=1}^{n}\Sign(E_{l}/G)\notag\\
 &=\sum_{l=1}^{n}\Sign(E_{l}/G)\label{basesp}
\end{align}
by the Novikov additivity. Substituting (\ref{G-sign}) and (\ref{basesp}) to (\ref{aver}), we see that the signature of $E$ localizes if and only if the right hand side of (\ref{G-sign}) localizes. By the definition, any vertical component $F$ is included in a singular fiber. Hence, the normal Euler number $\chi(N(F))$ depends only on the singular fiber germ. Thus, if the normal euler numbers of $I(h,\psi)$ localizes, the signature of $E$ also localizes.
\end{proof}
By the equations (\ref{aver}),(\ref{G-sign}), and (\ref{basesp}), we can write the local signature as Theorem \ref{thm:main}.
\begin{remark}
In the proof, we do not consider the localization of the normal euler number of each horizontal component but the sum of the normal euler numbers of the horizontal components $I(h,\psi)$. 
\end{remark}

\section{The local normal euler number of the set of horizontal components}\label{section:local euler}
Recall that $G$ is a finite group, and $p:\Sigma_g\to S^2$ is a $G$-covering. Denote by $m$ the order of the branch set of the $G$-covering $p$. In this section, we assume that the order $m$ is at least $3$. Let $f_0:E_0\to B_0$ be a $\Sigma_g$-bundle on a manifold $B_0$ with structure group in $C(p)$. The group $G$ acts on the total space $E_0$ fiberwise as stated after Definition \ref{broad fibration}. In Section \ref{sub:multi-section}, we review a canonical multi-section of the normal bundles of $E_0^h$ for $h\ne1\in G$ constructed by Furuta. As in Section \ref{section:local signature}, let $f:E\to B$ be a fibration of the $G$-covering $p$, and $\{f_l,E_l,\Delta_l\}_{l=1}^n$ representatives of the fiber germs in the fibration. Applying Section \ref{sub:multi-section} to the $\Sigma_g$-bundle $E-\amalg_{l=1}^n \Int E_l\to B-\amalg_{l=1}^n\Int\Delta_l$, we have the canonical multi-section of the normal bundle of $S\cap(E-\amalg_{l=1}^n \Int E_l)$ for each horizontal component $S$. We show that the normal euler number of $I(h,\psi)$ in the fibration $f:E\to B$ localizes in Section \ref{section:local euler number} by means of this multi-section. 

\subsection{Multi-sections of the normal bundles of the fixed point sets}\label{sub:multi-section}
We review a canonical multi-section of the normal bundle of the fixed point set $E_0^h$ for $h\ne1\in G$ constructed by Furuta. Let $J_0$ denote the union $\cup_{h\ne 1\in G}E_0^h$ of all fixed point sets. Denote by $\bar{E}_{0}$ and $\bar{J}_{0}$ the quotient spaces of $E_0$ and $J_0$ under the $G$-action, respectively. Note that the smooth map $\bar{f}_0:\bar{E}_0\to B_0$ induced by $f_0:E_0\to B_0$ is a $S^2$-bundle. Denote by $q:E_0\to\bar{E}_0$ the quotient map of the $G$-action. Fix a fiberwise complex structure on the $S^2$-bundle $\bar{f}_0:\bar{E}_0\to B_0$. Since the restriction of $q$ to each fiber $f^{-1}(b)\to \bar{f}^{-1}(b)$ is isomorphic to the $G$-covering $p$, it induces the fiberwise complex structure on the $\Sigma_g$-bundle $f_0:E_0\to B_0$. 

\begin{lemma}[Furuta \cite{furuta1999lsbse} Lemma 2]\label{lemma:multi-section}
When the order $m$ of the branch set in $S^2$ of the $G$-covering $p$ is at least $3$, there exists a canonical section $\bar{s}$ of the complex line bundle $T(\bar{E}_0/B_0)|_{\bar{J}_0}^{\otimes(m-1)(m-2)}\to \bar{J}_0$.
\end{lemma}

In fact, he also constructed a canonical multi-section of $T(\bar{E}_0/B_0)|_{\bar{E}_0-\bar{J}_0}$, which we do not need in this paper.
\begin{proof}
For $b\in B_{0}$, the intersection of $\bar{J}_0$ and the fiber $\bar{f}^{-1}(b)$ is a point set. Number it as $\{\alpha_{i}(b)\}_{i=1}^{m}=\bar{J}_0\cap\bar{f}^{-1}(b)$. We will construct a tangent vector at $\alpha_{i}(b)$ by choosing other two branched points $\alpha_j(b)$ and $\alpha_k(b)$. Define an isomorphism
\[
t_{b}^{ijk}:\mathbf{CP}^{1}\to\bar{f}^{-1}(b)
\]
by mapping $0$, $1$, and $\infty$, to $\alpha_{i}(b)$, $\alpha_{j}(b)$, and
$\alpha_{k}(b)$, respectively. Since $\bar{f}^{-1}(b)$ has the complex structure, this isomorphism is unique. In this way, we obtain the tangent vector ${t_{b}}_{*}^{ijk}(\frac{d}{dz})$ at $\alpha_{i}(b)$, where $z$ is the inhomogeneous coordinate in $\mathbf{CP}^{1}$. If we move $j,k\in\{1,2,\cdots,m\}$, we have 
\[
\bigotimes_{j,k}{t_{b}}_{*}^{ijk}\left(\frac{d}{dz}\right)\in T_{\alpha_{i}(b)}(\bar{E_{0}}/B_{0})^{\otimes(m-1)(m-2)}.
\]
Thus, we obtain the non-zero section $\bar{s}$ of the bundle $T(\bar{E}_{0}/B_{0})^{\otimes(m-1)(m-2)}|_{\bar{J}_{0}}$.
\end{proof}

By means of the section $\bar{s}$ above, Furuta constructed a multi-section of the normal bundle of $E_0^h$ for $h\ne1\in G$ as in the following theorem (See the proof of Theorem 1 in Furuta \cite{furuta1999lsbse}). For a connected component $S_0$ of the fixed point set $E^h$, denote by $r_{S_0}$ the order of the subgroup of $G$ which fixes $S_0$ in $E_0$ pointwise. 

\begin{theorem}[Furuta]\label{theorem:multi-section}
Let $S_0$ be a connected component of the fixed point set $E_0^h$ for $h\ne1\in G$. When the order $m$ of the branch set of the $G$-covering $p$ is at least $3$, there exists a canonical section $s:S_{0}\to T(E_{0}/B_0)^{\otimes r_{S_0}(m-1)(m-2)}|_{S_0}$. Moreover, the homotopy class of the section $s_{S_0}$ does not depend on the choice of the complex structure.
\end{theorem}

\begin{proof}
Let $\bar{S}_0$ denote the image of $S_0$ under the map $q$. The restriction $\bar{S}_{0}\to T(\bar{E}_{0}/B_{0})^{\otimes(m-1)(m-2)}|_{\bar{S}_0}$ of the section $\bar{s}$ induces $s'_{S_0}:S_{0}\to q^{*}T(\bar{E}_{0}/B_{0})^{\otimes(m-1)(m-2)}|_{S_0}$. The map $q$ induces the isomorphism
\[
L:T(E_{0}/B_{0})^{\otimes r_{S_0}}|_{S_0}\to q^{*}T(\bar{E}_{0}/B_{0})|_{S_0}.
\] 
Hence we have obtained the desired nonzero section $s=(L^{-1})^{\otimes(m-1)(m-2)}s'_{S_0}$ of the bundle $T(E_{0}/B_{0})^{\otimes r_{S_0}(m-1)(m-2)}|_{S_0}$. Denote it by $s_{S_0}$. Since the set of fiberwise complex structure on the $S^2$-bundle $\bar{f}_0$ is contractible, the homotopy class of $s$ does not depend on the choice of the complex structure. 
\end{proof}

Choose a $G$-equivariant Riemannian metric on $T(E_0/B_0)|_{J_0}$ which is compatible with the almost complex structure. Let $S_0$ be a connected component of $E^h$ for arbitrary $h\in G-\{1\}$. The Riemannian metric induces a metric on $T(E_0/B_0)^{\otimes r_{S_0}}|_{S_0}$. This bundle is isomorphic to the pullback of $T(\bar{E}_0/B_0)|_{\bar{S}_0}$ by the map $q|_{S_0}:S_0\to \bar{S}_0$. Hence we obtain the Riemannian metric on $T(\bar{E}_0/B_0)|_{\bar{J}_0}$ which is compatible with the almost complex structure. Let $g_{J_0}$ and $g_{\bar{J}_0}$ be arbitrary Riemannian metrics on $T{J_0}$ and $T\bar{J_0}$, respectively. The bundles $TE_0|_{J_0}$ and $T\bar{E}_0|_{\bar{J}_0}$ have splittings $TE_0|_{J_0}=T(E_0/B_0)|_{J_0}\oplus TJ_0$ and $T\bar{E}_0|_{\bar{J}_0}=T(\bar{E}_0/B_0)|_{\bar{J}_0}\oplus T\bar{J}_0$. We can choose Riemannian metrics $g_r$ on $E_0$ and $\bar{g}_r$ on $\bar{E}_0$ whose restrictions to $TE_0|_{J_0}$ and $T\bar{E}_0|_{\bar{J}_0}$ are equal to $g_{E_0/B_0}|_{J_0}\oplus g_{J_0}$ and $g_{\bar{E}_0/B_0}|_{\bar{J}_0}\oplus g_{\bar{J}_0}$, respectively. Since the bundles $TJ_0$ and $T(E_0/B_0)|_{J_0}$ are orthogonal, $T(E_0/B_0)|_{S_0}$ and $N(S_0)$ are canonically isomorphic. Thus, the section $s:S_0\to T(E_0/B_0)|_{S_0}^{\otimes r_{S_0}(m-1)(m-2)}$ induces a section of $N(S_0)^{\otimes r_{S_0}(m-1)(m-2)}$. Denote it as $s_{S_0}$. In the same way, by the isomorphism $T(\bar{E}_0/B_0)|_{\bar{S}_0}\cong N(\bar{S}_0)$, the section $\bar{s}|_{\bar{S}_0}:\bar{S}_0\to T(\bar{E}_0/B_0)|_{\bar{S}_0}^{\otimes(m-1)(m-2)}$ induces a section $s_{\bar{S}_0}$ of $N(\bar{S}_0)^{\otimes(m-1)(m-2)}$.

\subsection{The local euler number}\label{section:local euler number}
We prove that the normal euler number of the horizontal components $I(h,\psi)$ localizes. Let $S$ be a compact connected surface with nonempty boundary, and $V(S)\to S$ a vector bundle. Assume that $V(S)$ is oriented as a manifold, and that we are given a nonzero section $s:\partial S\to V(S)|_{\partial S}$. We introduce an integer $n(s,V(S))$ for the section $s$ in order to localize the euler number of $I(h,\psi)$. 

In the following, all homology groups are with integral coefficients if not specified. Let $s_0:S\to V(S)$ be the zero section. If we extend the section $s$ to the section $\tilde{s}$ of $V(S)\to S$, the exact sequence
\[
0=H_{2}(V(S))\to H_{2}(V(S),V(S)-s_0(S))\to H_{1}(V(S)-s_0(S))\to H_{1}(V(S))
\]
shows the homology class $[\tilde{s}]\in H_{2}(V(S),V(S)-s_0(S))$ is independent of the choice of $\tilde{s}$. Denote by $[s_{0}]$ the homology class of the zero section in $H_2(V(S),V(S)|_{\partial S})$. Let $D(S)$ be a unit disk bundle in $V(S)$, and $S(S)$ its sphere bundle. Then, we have natural isomorphisms $H_{2}(V(S),V(S)|_{\partial S})\cong H_2(D(S),D(S)|_{\partial S})$ and $H_{2}(V(S),V(S)-s_0(S))\cong H_{2}(D(S),S(S))$ induced by the inclusions. Define a number $n(s,V(S)):=[\tilde{s}]\cdot[s_{0}]\in\mathbf{Z}$ in terms of the intersection form $H_{2}(D(S),D(S)|_{\partial S})\times H_{2}(D(S),S(S))\to\mathbf{Z}$. 

Let $f:E\to \Delta$ be a representative of a fiber germ $[f,E,\Delta]\in S_g^p$. Endow a Riemannian structure $g_r$ on $E$ whose restriction on $TE|_{\partial J}$ is equal to $g_{\partial E/\partial \Delta}|_{\partial J}\oplus g_J|_{\partial J}$ for some Riemannian metric $g_J$ on $TJ$ as in the last paragraph of Section \ref{sub:multi-section}. By applying Theorem \ref{theorem:multi-section} to $\partial E\to\partial \Delta$, we have a canonical section $s_{\partial S}:\partial S\to N(S)^{\otimes r_S(m-1)(m-2)}|_{\partial S}$ for each horizontal component $S\subset I(h,\psi)$. The local euler number is described in terms of this section $s_{\partial S}$ as follows.

\begin{definition}\label{local euler number}
Define a map $\chi_{\loc}^{h,\psi}:S_{g}^{p}\to \mathbf{Q}$ by
\[
\begin{array}{cccc}
S_{g}^{p} & \to & \mathbf{Q}\\
{}[f,E,\Delta] & \mapsto & \displaystyle\sum_{S\subset I(h,\psi)}\frac{1}{r_S(m-1)(m-2)}n(s_{\partial S},N(S)^{\otimes r_S(m-1)(m-2)}),
\end{array}
\]
where $r_S$ is the order of the subgroup of $G$ which fixes $S$ in $E$ pointwise.
\end{definition}

This map is well-defined as shown in the following. We need to show that the number $\chi_{\loc}^{h,\psi}([f,E,B])\in\mathbf{Q}$ is independent of the choice of the metric $g_r$. If we choose another extension $g'_r$ on $TE$ of $g_0$ such that $g_r=g'_r$ on $TE|_{\partial E}$, the corresponding section $s_{\partial S}$ of the normal bundle $N(S)$ does not change. Since the space of Riemannian metrics on $TJ$ is contractible, the number $\chi_{\loc}^{h,\psi}([f,E,B])\in\mathbf{Q}$ is also independent of the choice of the Riemannian metric $g_{J}$. Hence the local euler number $\chi_{\loc}^{h,\psi}:S_{g}^{p}\to \mathbf{Q}$ is well-defined.

\begin{theorem}
For any fibration $f:E\to B$ of the $G$-covering $p$ with singular fiber germs $\{[f_{l},E_l,\Delta_l]\}_{l=1}^{n}$, we have
\[
\sum_{S\subset I(h,\psi)}\chi(N(S))=\sum_{l=1}^{n}\chi_{\loc}^{h,\psi}([f_{l},E_l,\Delta_l]).
\]
In short, the map $\chi_{\loc}^{h,\psi}:S_{g}^{p}\to \mathbf{Q}$ is a local normal euler number.
\end{theorem}

\begin{proof}
Let $\{f_l,E_l,\Delta_l\}_{l=1}^n$ be representatives of the fiber germs in the fibration $f:E\to B$. We may assume that $\Delta_l$ are mutually disjoint in $B$. Choose a Riemannian metric $g_0$ on $E-\amalg_{l=1}^n \Int E_l$ as in the last paragraph of Section \ref{sub:multi-section}. We can extend the Riemannian metric $g_0$ on $E-\amalg_{l=1}^n \Int E_l$ to a Riemannian metric $g_r$ on $TE$. Denote by $S_0$ the intersection of $S$ and the complement of $\amalg_{l=1}^n\Int E_l$. Let $\tilde{N}$ denote the tensor product $N(S)^{\otimes r_{S_0}(m-1)(m-2)}$ of the normal bundle of $S$ in $E$. By Theorem \ref{theorem:multi-section}, we have the canonical section $s_{S_0}$ of $\tilde{N}\to S$ over $S_0$. Extend the section $s_{S_0}$ to a section $\tilde{s}$ of $\tilde{N}\to S$ over $S$ transversely to the zero section $s_{0}:S\to \tilde{N}$. With the number $n(s_{S_0}|_{S_0 \cap\partial E_l},\tilde{N}|_{S\cap E_l})$, the intersection of $[\tilde{s}]\in H_{2}(\tilde{N},\tilde{N}-S)$ and $[s_{0}]\in H_{2}(\tilde{N})$ is described as
\[
\tilde{s}\cdot s_{0}=\sum_{l=1}^{n}n(s_{S_0}|_{S_0\cap\partial E_{l}},\tilde{N}|_{S\cap E_l})=\sum_{l=1}^{n}r_{S}(m-1)(m-2)\chi_{\loc}^{h,\psi}(f_{l}).
\]
Since the euler number of $\tilde{N}$ is equal to the self-intersection number of a section of $\tilde{N}\to S$, this is equal to $\chi(\tilde{N})=r_{S}(m-1)(m-2)\chi(N(S))$. Hence we have 
\[
\sum_{S\subset I(h,\psi)}\chi(N(S))=\sum_{l=1}^{n}\chi_{\loc}^{h,\psi}([f_{l},E_l,B_l]).
\]
\end{proof}

\section{The Meyer cocycle and symmetric mapping class groups}
Let $\Diff_{+}\Sigma_g$ denote the orientation-preserving diffeomorphism group of the closed surface $\Sigma_g$ of genus $g$. The mapping class group $\mathcal{M}_{g}$ of the surface $\Sigma_g$ is defined by the path-connected component $\pi_{0}\Diff_{+}\Sigma_{g}$ of this topological group with $C^{\infty}$ topology. For a finite regular covering $p:\Sigma_g\to\Sigma$ on a compact surface $\Sigma$, Birman-Hilden \cite{birman1973ihr} defined a group $\mathcal{M}_g(p)$ called the symmetric mapping class group. We restrict ourselves to the case $\Sigma=S^2$. Recall that we denote by $C(p)$ the centralizer of the deck transformation group of the $G$-covering $p$.

\begin{definition}
The symmetric mapping class group of the $G$-covering $p:\Sigma_g\to S^2$ is defined by
\[
\mathcal{M}_{g}(p):=\pi_{0}C(p).
\]
\end{definition}
Let $T$ be a finite set in $\Sigma_g$. Denote by $\Diff_{+}(\Sigma_{g},T)$ the group of orientation-preserving diffeomorphisms on the surface $\Sigma_g$ which fixes the set $T$ pointwise. Denote by $A=\{\alpha_{1},\alpha_{2},\cdots,\alpha_{m}\}\subset S^2$ the branch set of $p$. Pick a point $*$ in $S^2-A$. In the same way, the symmetric mapping class group of the pointed surface $(\Sigma_g,p^{-1}(*))$ is defined by $\mathcal{M}_{g}^{(*)}(p):=\pi_{0}(\Diff_{+}(\Sigma_{g},p^{-1}(*))\cap C(p))$.

Mapping a path-connected component of $C(p)$ to the corresponding component of $\Diff_{+}\Sigma_{g}$, we obtain the natural homomorphism
\[
\Phi:\mathcal{M}_{g}(p)\to\mathcal{M}_{g}.
\]
Assume that the $G$-covering $p:\Sigma_g\to S^2$ has at least $3$ branch points as in Section \ref{section:local euler}. Meyer \cite{meyer1973sf} introduced a 2-cocycle of the mapping class group $\mathcal{M}_g$, called the Meyer cocycle. We construct a cobounding function of the pullback of the Meyer cocycle by $\Phi$ in a subgroup $\mathcal{M}_{\mon}(p)\subset\mathcal{M}_g(p)$ in Theorem \ref{thm:singlar even}, when $p:\Sigma_g\to S^2$ has at least $3$ branch points. When this subgroup coincides with the symmetric mapping class group $\mathcal{M}_g(p)$, we can construct a local signature for fibrations of the $G$-covering $p$ in the broad sense, using this cobounding function (Proposition \ref{local signature broad}).

The symmetric mapping class group $\mathcal{M}_{g}(p)$ arises as the monodromy group of the $\Sigma_g$-bundles whose structure groups are included in $C(p)$. Let us recall the monodromy homomorphisms of $\Sigma_g$-bundles. Let $f:E\to B$ be a $\Sigma_g$-bundle on a manifold $B$. Fix a base point $b$ in $B$, and an identification $\Psi_0:\Sigma_g\to f^{-1}(b)$, which is called the reference fiber. For a homotopy class $\gamma\in\pi_1(B,b)$, choose a based loop $l:[0,1]\to B$ which represents $\gamma$. Since the pullback $q:l^*E\to [0,1]$ is the trivial $\Sigma_g$-bundle, we can choose a trivialization $\tilde{\Psi}:\Sigma_g\times [0,1]\to l^*E$ such that $\tilde{\Psi}(x,0)=\Psi_0(x)$. Define the diffeomorphism $\Psi_1:\Sigma_g\to f^{-1}(b)$ by $\Psi_1(x)=\tilde{\Psi}(x,1)$. The isotopy class of the diffeomorphism $\Psi_1^{-1}\Psi_0$ is called the monodromy of the $\Sigma_g$-bundle $f:E\to B$ along the loop $l$. It does not depend on the choice of $\tilde{\Psi}$ and $l$. Thus, we can define a homomorphism $\pi_1(B,b)\to \mathcal{M}_g$, called the monodromy homomorphism. If the structure group of the $\Sigma_g$-bundle $f:E\to B$ is in $C(p)$, the diffeomorphism $\Psi_1^{-1}\Psi_0$ is also included in $C(p)$. Similarly, we have the homomorphism $\pi_1(B,b)\to \mathcal{M}_g(p)$. We also call it the monodromy homomorphism.

Define the subgroup $\mathcal{M}_{\mon}(p)$ of the symmetric mapping class group $\mathcal{M}_g(p)$ as follows.

\begin{definition}
Denote by $\mathcal{M}_{\mon}(p)$ the subgroup of the symmetric mapping class group $\mathcal{M}_{g}(p)$ normally generated by the monodromies which arise in the set of fiber germs $S_{g}^{p}$ along the boundary circles.
\end{definition}

Let us review the definition of the Meyer cocycle. For $i=1,2,3$, let $D_{i}$ be disjoint closed disks in a $2$-sphere. Denote by $P$ a pair of pants $S^{2}-\amalg_{i=1}^{3}\Int D_{i}$. Let $E_g^{\varphi,\psi}$ be the total space of a $\Sigma_g$-bundle on $P$ whose monodromies around the boundary circles $\partial D_i$ are given by $\varphi,\psi,(\varphi\psi)^{-1}\in\mathcal{M}_g$. This $\Sigma_g$-bundle is unique up to isomorphism. 

\begin{definition}[Meyer \cite{meyer1973sf}]
The map
\[
\begin{array}{cccccc}
\tau_{g}: & \mathcal{M}_{g} & \times & \mathcal{M}_{g} & \to & \mathbf{Z},\\
 & (\varphi & , & \psi) & \mapsto & \Sign E_{g}^{\varphi,\psi}
 \end{array}
\]
is called the Meyer cocycle.
\end{definition}

Meyer proved that the map $\tau$ is a 2-cocycle on the mapping class group. Moreover, he showed that this cocycle represents a nontrivial 2-cohomology class of $\mathcal{M}_g$ when $g\ge3$. Later, it is rediscovered by Turaev \cite{turaev1987fsc}.

Birman-Hilden (Theorem 1 in \cite{birman1973ihr}) shows that if the deck transformation group fixes the branch set pointwise, then $\Phi$ is injective. Let $p':\Sigma_g\to S^2$ be a $\mathbf{Z}_2$-covering on $S^2$ for $g\ge2$. Especially, the symmetric mapping class group for $p'$ is isomorphic to a subgroup of the mapping class group called the hyperelliptic mapping class group.

To construct a cobounding function of $\Phi^*\tau_g$, we need a following lemma.

\begin{lemma}
For a mapping class $\hat{\varphi}\in\mathcal{M}_{\mon}(p)$, there exists a fibration $f:E\to D^2$ of the $G$-covering $p$ whose monodromy along the boundary circle $\partial D^2$ is $\hat{\varphi}$.
\end{lemma}

\begin{proof}
Let $f:E\to \Delta$ be a representative of a fiber germ. Let $t$ be a point in $\partial \Delta$, and $\Psi_0:\Sigma_g\to f^{-1}(t)$ a reference fiber. Assume that the monodromy of $f$ along the boundary circle is given by $\hat{\varphi}_0\in\mathcal{M}_{\mon}(p)$. First, we construct fiber germs whose monodromies are $\hat{\psi}\hat{\varphi}_0\hat{\psi}^{-1}$ and $\hat{\varphi}^{-1}_0$ for $\hat{\psi}\in\mathcal{M}_g(p)$.

Let $\hat{h}$ be a diffeomorphism of $\Sigma_g$ which represents $\hat{\varphi}\in\mathcal{M}_g(p)$. If we change the reference fiber by $\Psi_0 \hat{h}^{-1}:\Sigma_g\to f^{-1}(t)$, the monodromy is given by $\hat{\psi}\hat{\varphi}_0\hat{\psi}^{-1}$. Choose an orientation-reversing diffeomorphism $\iota:\Delta\to\Delta$. If we endow the other orientation on $E$, the smooth map $\iota f:E\to\Delta$ is also a fibration of the $G$-covering $p$. this fibration has $\hat{\varphi}^{-1}_0$ as its monodromy along the boundary circle. 

Let $f_i:E_i\to\Delta_i$ be representatives of fiber germs for $i=1,2,\cdots,n$. Choose a point $t_i$ in each boundary, and a reference fiber $\Psi_i:\Sigma_g\to f_i^{-1}(t_i)$. Assume that the monodromies of $f_i$ along the boundary circles are $\hat{\varphi}_i$ for $i=1,2,\cdots,n$. It suffices to construct a fibration whose monodromy along the boundary circle is given by $\prod_{i=1}^n \hat{\varphi}_i$. If we glue each reference fiber $f_i^{-1}(t_i)$ in $E_i$ together by $\Psi_{i'}\Psi_i^{-1}:f_i^{-1}(t_i)\to f_{i'}^{-1}(t_{i'})$ for $1\le i,i'\le n$, we obtain the space $\bigcup_{i=1}^n E_i$. Denote by $\bigvee_{i=1}^n \Delta_i$ the wedge sum obtained by identifying each $t_i\in\Delta_i$. The maps $f_i$ induce the map $F:\bigcup_{i=1}^n E_i\to \bigvee_{i=1}^n \Delta_i$. Embed the wedge sum $\bigvee_{i=1}^n \Delta_i$ in $D^2$ as in Figure \ref{fig:retract.eps}.
\begin{figure}[htbp]
  \begin{center}
    \includegraphics{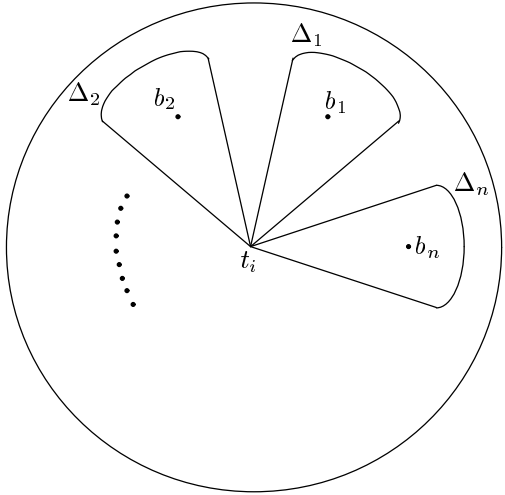}
  \end{center}
  \caption{the wedge sum $\bigvee_{i=1}^n \Delta_i$ in $D^2$}
  \label{fig:retract.eps}
\end{figure}
Then, there exists a deformation retraction $r:D^2\to\bigvee_{i=1}^n \Delta_i$. Denote by $b_i\in\Delta_i$ the critical value of $f_i$. By taking the pullback of the $\Sigma_g$-bundle $\bigcup_{i=1}^n(E_i-f_i^{-1}(b_i))\to \bigvee_{i=1}^n (\Delta_i-b_i)$ by $r$, we obtain the fibration on $D^2$ whose monodromy along the boundary circle is $\prod_{i=1}^n \hat{\varphi}_i$. There exists a (smooth) fibration topologically isomorphic to this fibration.
\end{proof}

\begin{theorem}\label{thm:singlar even}
Let $G$ be a finite group, and $p:\Sigma_g\to S^2$ a $G$-covering with at least $3$ branch points. For a mapping class $\hat{\varphi}\in\mathcal{M}_{\mon}(p)$, there exists a fibration $f:X\to D^2$ of the $G$-covering $p$ on a closed 2-disk whose monodromy along the boundary circle $\partial D^2$ is $\hat{\varphi}$. Denote by $\{[f_{l}, E_l, \Delta_l]\}_{l=1}^n$ the fiber germs arise in the fibration. Define a map $\phi:\mathcal{M}_{\mon}(p)\to\mathbf{Q}$ by
\[
\phi(\hat{\varphi}):=\sum_{l=1}^{n}\sigma_{\loc}([f_{l},E_l,\Delta_l])-\Sign X,
\]
where $\sigma_{\loc}$ is the local signature described in Theorem \ref{thm:main}. This is well-defined, and cobounds the 2-cocycle $\Phi^{*}\tau_{g}$ in $\mathcal{M}_{\mon}(p)$. That is to say, it satisfies
\[
\phi(\hat{\varphi})+\phi(\hat{\psi})+\phi((\hat{\varphi}\hat{\psi})^{-1})=\Phi^{*}\tau_{g}(\hat{\varphi},\hat{\psi}),
\]
for $\hat{\varphi},\hat{\psi}\in\mathcal{M}_{\mon}(p)$.
\end{theorem}

\begin{proof}
Choose a $\Sigma_g$-bundle $E\to P$ on the pair of pants whose monodromies along the boundary circles are $\hat{\varphi}$, $\hat{\psi}$, and $(\hat{\varphi}\hat{\psi})^{-1}$, respectively. Since these mapping classes lie in $\mathcal{M}_{\mon}(p)$, there exist fibrations $X_{i}\to D^{2}$ ($i=1,2,3$) of the $G$-covering $p$ whose monodromies along the boundary circles are $\hat{\varphi}$, $\hat{\psi}$, and $(\hat{\varphi}\hat{\psi})^{-1}$, respectively. Let $\{f_{l}^{i}, E_l^i, \Delta_l^i\}_{l=1}^{n_{i}}$ denote representatives of the fiber germs arise in the fibration $X_i$. By the definition of the Meyer cocycle, we have
\[
\Phi^{*}\tau_{g}(\hat{\varphi},\hat{\psi})=\Sign E.
\]
By the definition of the local signature, we have
\[
\Sign E=\sum_{i=1}^{3}\left(\sum_{l=1}^{n_{i}}\sigma_{\loc}([f_{l}^{i},E_l^i,\Delta_l^i])-\Sign X_{i}\right)=\phi(\hat{\varphi})+\phi(\hat{\psi})+\phi((\hat{\varphi}\hat{\psi})^{-1}).
\]
If we substitute another fibration $X_1'\to D^2$ whose monodromy along the boundary circle is $\hat{\varphi}$ for $X_1\to D^2$, the left-hand side of the equation does not change. Hence the value $\phi(\hat{\varphi})$ does not also depend on the choice of $X_1$. Moreover, the equation shows that the map $\phi$ cobounds the Meyer cocycle.
\end{proof}

In terms of the cobounding function, we obtain a local signature for fibrations of the $G$-covering $p$ in the broad sense, if the subgroup $\mathcal{M}_{\mon}(p)$ coincides with the whole group $\mathcal{M}_g(p)$. This local signature is defined as a function on another kind of fiber germs. Denote by $\Delta$ a closed 2-disk. Consider fibrations ($f$,$E$,$\Delta$) of the $G$-covering $p$ in the broad sense with unique critical values $b\in\Delta$. Let ($f_1$,$E_1$,$\Delta_1$) and ($f_2$,$E_2$,$\Delta_2$) be such fibrations which have unique critical values $b_1$ and $b_2$, respectively. These fibrations are equivalent if and only if there exist
\begin{enumerate}
\item closed 2-disks $\Delta_1'\subset\Delta_1$ and $\Delta_2'\subset\Delta_2$ including the critical values,
\item an orientation-preserving diffeomorphism $\varphi:(\Delta_1',b_1)\to(\Delta_2',b_2)$,
\item an orientation-preserving diffeomorphism $\tilde{\varphi}:f_1^{-1}(\Delta_1)\to f_2^{-1}(\Delta_2)$ such that $\varphi f_1=f_2\tilde{\varphi}$ and it restricts to a $G$-equivariant diffeomorphism $f_1^{-1}(\Delta_1-b_1)\to f_2^{-1}(\Delta_2-b_2)$.
\end{enumerate}
We denote this set of equivalent classes by $\tilde{S}_g^p$. 

\begin{proposition}\label{local signature broad}
Let $G$ be a finite group, and $p:\Sigma_g\to S^2$ a $G$-covering with at least 3 branch points. Assume that the group $\mathcal{M}_{\mon}(p)$ coincides with the whole group $\mathcal{M}_g(p)$. Let $\phi:\mathcal{M}_g(p)\to\mathbf{Q}$ be the cobounding function of the pullback $\Phi^{*}\tau_{g}$ of the Meyer cocycle in Theorem \ref{thm:singlar even}. For a fiber germ $[f,E,\Delta]\in\tilde{S}_g^p$, denote by $\hat{\varphi}$ the monodromy along the boundary curve $\partial\Delta$. The map $\sigma'_{\loc}:\tilde{S}_g^p\to\mathbf{Q}$ defined by
\[
\sigma'_{\loc}([f,E,\Delta]):=\phi(\hat{\varphi})+\Sign E
\]
is a local signature for fibrations of the $G$-covering $p$ in the broad sense. 
\end{proposition}
The proof is the same as that of Theorem 4.4 in Endo \cite{endo2000mss}.

\begin{remark}
In general, a cobounding function of the pullback $\Phi^{*}\tau_{g}$ of the Meyer cocycle in $\mathcal{M}_{\mon}(p)$ is not unique.
\end{remark}

\section{Generators of symmetric mapping class groups}\label{section:generator}
In this section, we describe a generating set of the symmetric mapping class group $\mathcal{M}_g(p)$ of the $G$-covering $p$, assuming that the finite group $G$ is abelian. Let $d\ge2$ be an integer. In Section \ref{sec:fiber germs}, we will construct fiber germs whose monodromies are inverses of the standard generator system $\{\sigma_{ij}\}$, and prove $\mathcal{M}_g(p_1)=\mathcal{M}_{\mon}(p_1)$.

Let $\hat{f}$ be in the centralizer $C(p)$. The diffeomorphism $\hat{f}$ induces a diffeomorphism $f$ of $S^2$ satisfying a commutative diagram
\[
\begin{CD}\Sigma_{g}@>\hat{f}>>\Sigma_{g}\\
@VpVV@VpVV\\
S^{2}@>f>>S^{2}.
\end{CD}
\]
We call the diffeomorphism $f$ the projection of $\hat{f}$. Recall that $A$ is the branch set in $S^2$ of the $G$-covering $p$. Let $k:\pi_{1}(S^{2}-A)\to G$ be the monodromy homomorphism of the $G$-covering $p$. Since $G$ is abelian, this induces the homomorphism $\bar{k}:H_1(S^2-A)\to G$. Choose a base point $*$ in $S^2-A$. Denote by $\gamma_{\alpha_{i}}:[0,1]\to S^2$ a based loop which rotates around a point $\alpha_{i}$ counterclockwise once. For $h\in G-\{1\}$, define a subset $A_h$ of the branch set $A$ in $S^2$ by $A_h=\{\alpha\in A\,|\, \bar{k}_*[\gamma_{\alpha}]=h\}$, where $[\gamma_{\alpha}]$ is the homology class of $\gamma_{\alpha}$. Let $\mathcal{M}_{0}^{A}$ denote the mapping class group which preserves the set $A_h$ setwise for any $h\in G-\{1\}$. We also denote by $\mathcal{M}_{0}^{A,*}$ the mapping class group which preserves the base point $*$ and each $A_h$ setwise. The projection $f$ of $\hat{f}\in C(p)$ preserves each branch set $A_h$. For the details, see Proposition 1.2 in \cite{sato2009asm}. Thus, we have homomorphisms $\Phi':\mathcal{M}_{g}^{(*)}(p)\to \mathcal{M}_{0}^{A,*}$ and $\Phi:\mathcal{M}_{g}(p)\to \mathcal{M}_{0}^{A}$ defined by $[\hat{f}] \mapsto [f]$.

\begin{lemma}\label{lem: surj}
Assume that the finite group $G$ is abelian.
\begin{enumerate}
\item The homomorphism $\Phi':\mathcal{M}_{g}^{(*)}(p)\to \mathcal{M}_{0}^{A,*}$ is isomorphic. 
\item The homomorphism $\Phi:\mathcal{M}_{g}(p)\to \mathcal{M}_{0}^A$ is surjective, and the kernel is generated by $\Deck(p)$.
\end{enumerate}
\end{lemma}

\begin{figure}[htbp]
  \begin{center}
    \includegraphics{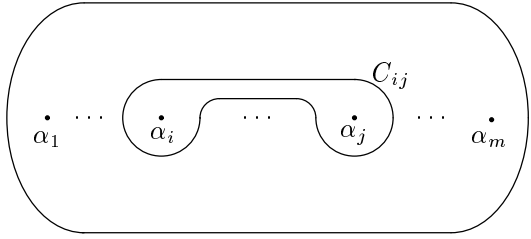}
  \end{center}
\caption{the simple closed curve $C_{ij}$}
\label{fig:sigma-ij.eps} 
\end{figure}

\begin{proof}
\if01
First we prove that the homomorphisms $\Phi$ and $\Phi'$ are surjective. Denote by $\rho:\pi_{1}(S^{2}-A)\to H_{1}(S^{2}-A;\mathbf{Z})$ the Hurewicz homomorphism. Since $G$ is abelian, the monodromy homomorphism $k:\pi_{1}(S^{2}-A,*)\to G$ factors through $H_{1}(S^{2}-A;\mathbf{Z})$. Thus, we have
\[
\Ker(\pi_{1}(S^{2}-A)\to G)=\rho^{-1}\Ker(H_{1}(S^{2}-A;\mathbf{Z})\to G).
\]
Since any mapping class $[f]$ in $\mathcal{M}_0^{A,*}$ (resp. $[f]$ in $\mathcal{M}_0^{A}$) acts on $H_{1}(S^{2}-A;\mathbf{Z})$ trivially, the lifting property of covering spaces implies that the diffeomorphism $f$ induces a diffeomorphism on the covering space $\Sigma_g$. Hence $\Phi$ and $\Phi'$ are surjective. 
\fi
The surjectivity of $\Phi$ is a special case of Proposition 1.2 in \cite{sato2009asm}. We can show that the homomorphism $\Phi'$ is surjective in the same way.

We compute the kernels of $\Phi$ and $\Phi'$. Let $f$ be the projection of $\hat{f}\in C(p)$. If the mapping class $[\hat{f}]\in \mathcal{M}_g(p)$ is in the kernel of $\Phi$, there exists an isotopy $\{f_{s}\}_{0\le s\le 1}$ satisfying $f_0=f$ and $f_1=id$. Choose the lift of this isotopy $\{\hat{f}_{s}\}_{0\le s\le 1}$ such that $\hat{f}_0=\hat{f}$. Since $\hat{f}_1$ is a lift of the identity map, it is a deck transformation. Hence, the kernel is generated by $\Deck(p)$. By the same argument, we can show that the kernel of $\Phi'$ is also generated by $\Deck(p)\cap \Diff_+(\Sigma_g,p^{-1}(*))$, which is the trivial group.
\end{proof}
For mutually distinct integers $i,j$, choose a simple closed curve $C_{ij}$ as in Figure \ref{fig:sigma-ij.eps}. When there exist mutually distinct $h,h'\in G-\{1\}$ which satisfy $\alpha_i\in A_{h}$ and $\alpha_j\in A_{h'}$, Denote by $\tau_{ij}\in\mathcal{M}_{0}^{A,*}$ the full Dehn twist along $C_{ij}$. When there exists $h\in G-\{1\}$ which satisfies $\alpha_{i},\alpha_{j}\in A_{h}$, denote by $\sigma_{ij}\in\mathcal{M}_{0}^{A,*}$ the half Dehn twist along $C_{ij}$. This is the mapping class which exchanges the points $\alpha_{i}$ and $\alpha_{j}$ and whose square is the full Dehn twist along $C_{ij}$. 
We denote by $\hat{\sigma}_{ij}$ and $\hat{\tau}_{ij}$ in $\mathcal{M}_{g}^{(*)}(p)$ the lifts of $\sigma_{ij}$ and $\tau_{ij}$ which preserve the fiber $p^{-1}(*)$ pointwise, respectively. The inclusion $\Diff_{+}(\Sigma_{g},p^{-1}(*))\cap C(p)\to C(p)$ induces a homomorphism $\mathcal{M}_{g}^{(*)}(p)\to\mathcal{M}_{g}(p)$. We also denote by the same symbol the images of $\hat{\sigma}_{ij}$ and $\hat{\tau}_{ij}$ under this homomorphism.

\begin{lemma}\label{lem:generator}
If the finite group $G$ is abelian, then both of the groups $\mathcal{M}_{g}^{(*)}(p)$ and $\mathcal{M}_{g}(p)$ are generated by $\hat{\sigma}_{ij}$ and $\hat{\tau}_{ij}$ for $1\le i\le m$ and $1\le j\le m$.
\end{lemma}

\begin{proof}
Since we have the isomorphism $\Phi':\mathcal{M}_{g}^{(*)}(p)\cong\mathcal{M}_{0}^{A,*}$, it suffices to show that
\begin{enumerate}
\item the mapping classes $\sigma_{ij}$, $\tau_{ij}$ generates $\mathcal{M}_{0}^{A,*}$,
\item the homomorphism $\mathcal{M}_{g}^{(*)}(p)\to\mathcal{M}_{g}(p)$ is surjective.
\end{enumerate}
First, we show (i). Let $\Diff_{+}(S^{2},A,*)$ be the diffeomorphism group which preserves the base point $*$ and the set $A$ pointwise. Denote by $n(h)$ the order of $A_{h}$, and by $S_{n(h)}$ the symmetric group of degree $n(h)$. Since $\mathcal{M}_0^{A,*}$ permutes the elements of each set $A_h$, we have the homomorphism $\eta:\mathcal{M}_{0}^{A,*}\to\prod_{h\in G-\{1\}}S_{n(h)}$. Since the group $\prod_{h\in G-\{1\}}S_{n(h)}$ is generated by the images of $\hat{\sigma}_{ij}$ under $\eta$, we have the exact sequence
\[
\begin{CD}
1@>>>\pi_{0}\Diff_{+}(S^{2},A,*)@>>>\mathcal{M}_{0}^{A,*}@>\eta>>\prod_{h\in G-\{1\}}S_{n(h)}@>>>1.
\end{CD}
\]
It is known that $\pi_{0}\Diff_{+}(S^{2},A,*)$ is generated by $\sigma_{ij}^2$ and $\tau_{ij}$. For example, see Section 1.5 in Birman \cite{birman1975bla}. Thus, we have proved (i).

Next, we show (ii). Let $\hat{*}$ be a point in $p^{-1}(*)$. Let $\hat{f}\in C(p)$ be a diffeomorphism, and let $f\in \Diff_+ S^2$ denote the projection of $\hat{f}$. The map $\rho:\Diff_+(S^2, A)\to S^2-A$ defined by $h\mapsto h(*)$ is a fiber bundle with fiber $\Diff_+(S^2,A,*)$ as in Theorem 4.1 in Birman \cite{birman1975bla}. Pick a path $\hat{\gamma}:[0,1]\to \Sigma_g-p^{-1}(A)$ such that $\hat{\gamma}(0)=\hat{*}$ and $\hat{\gamma}(1)=\hat{f}(\hat{*})$. Denote by $\Psi:[0,1]\to\Diff_+(S^2, A)$ a lift of $p\hat{\gamma}:[0,1]\to S^2-A$ with respect to the fiber bundle $\rho:\Diff_+(S^2, A)\to S^2-A$ such that $\Psi(0)$ is the identity map. By the lifting property of the $G$-covering $p:\Sigma_g\to S^2$, this can be lifted to the map $\hat{\Psi}:[0,1]\to C(p)$ such that $\hat{\Psi}(0)$ is the identity map. Then, the composite $\hat{\Psi}(1)^{-1}\hat{f}$ preserves the point $\hat{*}$. Moreover, since $\hat{\Psi}(1)^{-1}\hat{f}$ commutes with the deck transformations, it preserves each point of $p^{-1}(*)$. Hence $\hat{\Psi}(1)^{-1}\hat{f}$ represents a mapping class in $\mathcal{M}_g^{(*)}(p)$. By the isotopy $\hat{\Psi}$, we have $[\hat{\Psi}(1)^{-1}\hat{f}]=[\hat{f}]\in\mathcal{M}_{g}(p)$. This shows that $\mathcal{M}_{g}^{(*)}(p)\to\mathcal{M}_{g}(p)$ is surjective.
\end{proof}

\section{The construction of fiber germs}\label{sec:fiber germs}
Let $d\ge2$ and $m\ge3$ be integers such that $m$ is divided by $d$, and $A$ a finite set $\{\alpha_i\}_{i=1}^m$ in $S^2$. For each $i=1,2,\cdots,m$, choose a loop $\gamma_{\alpha_i}$ which rotates around a point $\alpha_i$ counterclockwise once. Define a surjective homomorphism $k:H_{1}(S^{2}-A)\to\mathbf{Z}_{d}$ by mapping each homology class $[\gamma_{\alpha_i}]$ to $1 \mod d$. Since $m\equiv0\mod d$, this is well-defined. Let $p_1:\Sigma_g\to S^2$ be the $\mathbf{Z}_d$-covering on $S^2$ which has the branch set $A$ in $S^{2}$ and the monodromy homomorphism $k$. In this section, we will construct a fiber germ of $\mathbf{Z}_d$-covering $p_1$ whose monodromy is the inverse$\hat{\sigma}_{12}^{-1}$ of the generator of $\mathcal{M}_g(p_1)$ intruduced in Section \ref{section:generator}. We will also calculate the local signature of this fiber germ and the value $\phi(\hat{\sigma}_{12})$ of the cobounding function of the pullback of the Meyer cocycle $\tau_g$ by the homomorphism $\Phi:\mathcal{M}_g(p)\to\mathcal{M}_g$.

\begin{lemma}\label{cover p1}
The subgroup $\mathcal{M}_{\mon}(p_1)$ coincides with the whole symmetric mapping class group $\mathcal{M}_g(p_1)$.
\end{lemma}

\begin{proof}
We need to show that the generating set $\{\hat{\sigma}_{ij}\}_{i,j\in A}$ of the symmetric mapping class group $\mathcal{M}_g(p_1)$ is contained in the subgroup $\mathcal{M}_{\mon}(p_1)$. Since $\{\hat{\sigma}_{ij}\}_{i,j\in A}$ are in the same conjugacy class in $\mathcal{M}_g(p_1)$, it suffices to construct a fiber germ in $S_g^{p_1}$ whose monodromy is $\hat{\sigma}_{12}$. 
\if01
Let $[x:y:z]$ be a homogeneous coordinate of $\mathbf{CP}^{2}$. Consider the smooth surface
\[
E=\{b_2y_1^{d}=y_2^d(x_1^{m'}-x_2^{m'})(b_2x_1^2-b_1x_2^2)\,|\,([x_1:x_2],[y_1:y_2],[b_1,b_2])\in\mathbf{CP}^1\times\mathbf{CP}^1\times\mathbf{CP}^{1}\}.
\]
Let $t$ be a generator $1$ in $\mathbf{Z}_{d}$. The surface $E$ has a $\mathbf{Z}_d$-action defined by
\[
t:([x_1:x_2],[y_1:y_2],[b_1,b_2])\mapsto([x_1:x_2],[e^{2\pi\sqrt{-1}/d} y_1:y_2],[b_1,b_2]).
\]
Define a smooth map $f:E\to\mathbf{CP}^1$ by $([x_1:x_2],[y_1:y_2],[b_1,b_2])\mapsto [b_1,b_2]$. Then, $f:E\to\mathbf{CP}^1$ is a fibration of the $\mathbf{Z}_d$-covering $p_1$. 
\fi

Let $\bar{E}$ be the product space $\Delta_1\times\mathbf{CP}^{1}$, where $\Delta_1=\{b\in \mathbf{C}\,|\,|b|\le 1/2\}$ is a closed 2-disk. Let $[x_{1}:x_{2}]$ be the homogeneous coordinate of $\mathbf{CP}^1$, and $m'=m-2$. Denote by $\bar{J}$ a submanifold of $\bar{E}$ defined by the equation $(x_{1}^{m'}-x_{2}^{m'})(x_{1}^{2}-bx_{2}^{2})=0$. For $b=re^{\sqrt{-1}\phi}\in\Delta_1$ where $r\ge0$ and $0\le\phi<2\pi$, define the root by $\sqrt{b}=\sqrt{r}e^{\sqrt{-1}\phi/2}$. Pick a diffeomorphism $T:S^2\to\{1/2\}\times\mathbf{CP}^{1}$ which maps $\alpha_1$, $\alpha_2$, and $\alpha_i$, to $(1/2,[\sqrt{1/2}:1])$, $(1/2,[-\sqrt{1/2}:1])$, and $(1/2,[1:e^{2\pi\sqrt{-1}(i-2)/m'}])$, for $i=3,4,\cdots,m$. Then, we have
\[
H_{1}(\bar{E}-\bar{J})=\bigoplus_{j=1}^m \mathbf{Z} e_j/(\mathbf{Z}(e_{1}-e_{2})\oplus\mathbf{Z}(e_{1}+e_{2}+\cdots+e_{m})),
\]
where $e_j$ is the homology class represented by the loop $T\circ\gamma_{\alpha_j}$. Define a homomorphism $l:H_{1}(\bar{E}-\bar{J})\to \mathbf{Z}_{d}$ by mapping each $e_j$ to $1$. This is well-defined since $m\equiv0\mod d$. Hence there exists a $\mathbf{Z}_{d}$-covering $q:E\to\bar{E}$ branched along $\bar{J}$ whose monodromy homomorphism is $l$. 

Denote by $\bar{f}:\bar{E}\to\Delta_1$ the projection to the second factor, and by $f:E\to\Delta_1$ the composite of $q:E\to\bar{E}$ and $\bar{f}$. The map $f$ is a fibration of $\mathbf{Z}_d$-covering $p_1$ in the narrow sense with the unique singular value $b=0$. We only have to check that it has the monodromy $\hat{\sigma}_{12}$. The restriction of $\bar{f}$
\[
(\bar{E}-(\{0\}\times\mathbf{CP}^1), \bar{J}-\{(0,[0:1])\})\to\Delta_1-0
\]
can be considered as a $m$-pointed sphere bundle. Consider $T:S^2\to\{1/2\}\times\mathbf{CP}^1$ as a reference fiber of this bundle. Then, the monodromy of this bundle along $\partial \Delta_1$ is $\sigma_{12}^{-1}\in\mathcal{M}_0^{A}$. Moreover, since there is a section $s:\Delta_1\to \bar{E}$ defined by $s(b)=(b,[1:0])$, the monodromy can be considered as $\sigma_{12}\in\mathcal{M}_0^{A,*}$. Since the isomorphism $\Phi':\mathcal{M}_g^{(*)}(p)\cong \mathcal{M}_0^{A,*}$ maps $\hat{\sigma}_{12}$ to $\sigma_{12}$, the fibration $f:E\to\Delta_1$ has the monodromy $\hat{\sigma}_{12}^{-1}\in\mathcal{M}_g(p)$. 
\end{proof}

\begin{lemma}
In the symmetric mapping class group $\mathcal{M}_g(p_1)$, the cobounding function of the pullback $\Phi^*\tau_g$ of the Meyer cocycle under $\Phi:\mathcal{M}_g(p_1)\to\mathcal{M}_g$ is unique.
\end{lemma}

\begin{proof}
If there exist two cobounding functions $\phi$ and $\phi'$, the map $\phi-\phi':\mathcal{M}_{g}(p_1)\to \mathbf{Q}$ is a $1$-cocycle. Hence it suffices to show that $H^1(\mathcal{M}_g(p_1);\mathbf{Q})=0$. By Lemma \ref{lem: surj}, we have the exact sequence
\[
\begin{CD}
\Deck(p_1)@>>>\mathcal{M}_{g}(p_1)@>>>\mathcal{M}_{0}^{A}@>>>1.
\end{CD}
\]
Since the deck transformation group $\Deck(p_1)$ is finite, we have the isomorphism
\[
H^1(\mathcal{M}_{0}^{m};\mathbf{Q})\cong H^1(\mathcal{M}_{g}(p_1);\mathbf{Q}).
\]
It is known that $H^1(\mathcal{M}_{0}^{m};\mathbf{Q})=0$, for example, this follows from the presentation of $\mathcal{M}_{0}^{m}$ obtained by Birman (Theorem 4.5 in \cite{birman1975bla}). Hence we have $H^1(\mathcal{M}_g(p_1);\mathbf{Q})=0$.
\end{proof}

By computing the local euler number of the fixed point set in $E$ of the $\mathbf{Z}_d$-action, we can calculate the cobounding function of the pullback of the Meyer cocycle. As in Proposition \ref{proposition:localsignature}, if $d=2$ i.e. $m=2g+2$, the value $\phi(\hat{\sigma}_{ij})$ for $\hat{\sigma}_{ij}\in\mathcal{M}_{g}(p_1)$ coincides with that of the Meyer function on the hyperelliptic mapping class group obtained by Endo (\cite{endo2000mss} Lemma 3.2). 

To prove Proposition \ref{proposition:localsignature}, we need Lemma \ref{intersection}. Let $D(\Delta)\to \Delta$ be a $D^2$-bundle on a closed $2$-disk $\Delta$, and $S(\Delta)$ its sphere bundle. The manifold $S(\Delta)$ induces the orientation on the boundary $S(\Delta)|_{\partial\Delta}$. Let $s:\Delta\to S(\Delta)$ be a section. For a section $s':\partial\Delta\to S(\Delta)|_{\partial \Delta}$, take an extension $\tilde{s'}:\Delta\to D(\Delta)$. Then we can consider the following two intersection numbers $H_{2}(D(\Delta),S(\Delta))\times H_{2}(D(\Delta),D(\Delta)|_{\partial \Delta})\to\mathbf{Z}$ and $H_{1}(S(\Delta)|_{\partial \Delta})\times H_{1}(S(\Delta)|_{\partial \Delta})\to\mathbf{Z}$. We denote by $\tilde{s'}\cdot s$ the former one, by $s'\cdot s|_{\partial \Delta}$ the latter one. Then

\begin{lemma}\label{intersection}
\[
\tilde{s'}\cdot s=-s'\cdot s|_{\partial \Delta}.
\]
\end{lemma}

\begin{proof}
We may identify the disk $\Delta$ with the embedded 2-disk $D^2=\{b\in \mathbf{C}\,|\, |b|\le 1\}$ in $\mathbf{C}$. The section $s$ gives a trivialization $D(\Delta)\cong D^{2}\times D^{2}$.

For some integer $k$, the section $s'$ represents the same class as the curve $\{(z,z^{k})|z\in S^{1}\}$ in $D^2\times D^2$ in $H_{1}(S(D^{2})|_{\partial D^{2}})$. Then, we have $s'\cdot s|_{\partial D^{2}}=-k$. Since the homomorphism $H_{2}(D(D^{2}),S(D^{2}))\to H_{1}(S(D^{2})|_{\partial D^{2}})$ is injective, the homology class in $H_{2}(D(D^{2}),S(D^{2}))$ of $\tilde{s'}$ is represented by the surface $\{(z,z^{k})|z\in D^{2}\}$. Hence we have
$\tilde{s'}\cdot s=k=-s'\cdot s|_{\partial D^{2}}$. 
\end{proof}

\begin{proof}[proof of Proposition \ref{proposition:localsignature}]
Let $E\to \Delta_1$ denote the fibration of the $\mathbf{Z}_d$-covering $p_1$ in Lemma \ref{cover p1}. Let $\bar{S}_{12}$ and $\bar{S}_{i}$ be the submanifolds defined by $x_{1}^{2}=bx_{2}^{2}$ and by $e^{2\pi\sqrt{-1}(i-2)/m'}x_{1}=x_{2}$ for $i=3,\cdots,m$, respectively. Denote by $S_{12}$ and $\{S_i\}_{i=3}^m$ their inverse images $q^{-1}(\bar{S}_{12})$ and $\{q^{-1}(\bar{S}_{i})\}_{i=3}^m$ under the map $q: E\to\bar{E}$. 

Let $r_{S_i}:TE|_{S_i}\to N(S_i)$, $r_{S_{12}}:TE|_{S_{12}}\to N(S_{12})$, $r_{\bar{S}_i}:T\bar{E}|_{\bar{S}_i}\to N(\bar{S}_i)$, and $r_{\bar{S}_{12}}:T\bar{E}|_{\bar{S}_{12}}\to N(\bar{S}_{12})$ be the natural projections. As in the last paragraph of Section \ref{sub:multi-section}, these projections induces the canonical isomorphisms
\begin{gather*}
N(S_i)|_{\partial S_i}\cong T(\partial E/\partial\Delta_1)|_{\partial S_i}, N(S_{12})|_{\partial S_{12}}\cong T(\partial E/\partial\Delta_1)|_{\partial S_{12}},\\
N(\bar{S}_i)|_{\partial\bar{S}_i}\cong T(\partial \bar{E}/\partial\Delta_1)|_{\partial\bar{S}_i}, \text{and } N(\bar{S}_{12})|_{\partial\bar{S}_{12}}\cong T(\partial \bar{E}/\partial\Delta_1)|_{\partial\bar{S}_{12}}.
\end{gather*}
Endow complex structures on $N(S_i)$, $N(S_{12})$, $N(\bar{S}_i)$, and $N(\bar{S}_{12})$ which are compatible with the inner products of the normal bundles.

In Section \ref{sub:multi-section}, we constructed the sections of $N(S_i)|_{\partial S_i}^{\otimes d(m-1)(m-2)}$, $N(S_{12})|_{\partial S_{12}}^{\otimes d(m-1)(m-2)}$, $N(\bar{S}_i)|_{\partial \bar{S}_i}^{\otimes(m-1)(m-2)}$, and $N(\bar{S}_{12})|_{\partial\bar{S}_{12}}^{\otimes(m-1)(m-2)}$, named $s_{\partial S_{12}}$, $s_{\partial S_i}$, $s_{\partial \bar{S}_{12}}$ and $s_{\partial \bar{S}_i}$, respectively. We review the definitions of the sections $s_{\partial \bar{S}_{12}}$ and $s_{\partial \bar{S}_i}$. Define a map $\alpha_i: \Delta_1\to \bar{E}$ by $\alpha_{1}(b)=(b,[1:\sqrt{b}])$, $\alpha_{2}(b)=(b,[1:-\sqrt{b}])$, and $\alpha_{i}(b)=(b,[1,e^{2\pi\sqrt{-1}(i-2)/m'}])$ for $i=3,\cdots,m$. Let $j$ and $k$ be integers such that $1\le j\le m$, $1\le k\le m$, and $i,j,k$ are mutually distinct. For such $j$ and $k$, define a not necessarily continuous section $s_i(j,k): \partial\bar{S}_i\to T(\bar{E}/\Delta_1)|_{\partial \bar{S}_i}$ by
\[
s_i(j,k)(\alpha_i(b))=(t_{b}^{ijk})_*\left(\frac{d}{dz}\right)
\]
as in Section \ref{sub:multi-section}. Note that, if $\{j,k\}\cap\{1,2\}\ne \emptyset$, the section $s_i(j,k)$ is not continuous since the root of $b$ is not continuous. The (continuous) section $s_{\partial \bar{S}_i}$ is defined by 
\[
s_{\partial \bar{S}_i}=r_{\bar{S}_i}(\bigotimes_{j,k}s_i(j,k)),
\]
where $j$ and $k$ run through integers such that $1\le j\le m$, $1\le k\le m$, and $i$, $j$, $k$ are distinct. In the same way, for integers $j$ and $k$ such that $3\le j\le m$, $3\le k\le m$, and $i,j,k$ are distinct, define sections of the bundle $T(\bar{E}/\Delta_1)|_{\partial \bar{S}_{12}}\to\bar{S}_{12}$ by 
\begin{gather*}
s_{12}(j,k)(\alpha_i(b))=(t_{b}^{ijk})_*\left(\frac{d}{dz}\right),\text{ for }i=1,2,\\
s_{12}^+(j)(\alpha_i(b))=
\begin{cases}
(t_{b}^{12j})_*\left(\frac{d}{dz}\right), &\text{ if }i=1,\\
(t_{b}^{21j})_*\left(\frac{d}{dz}\right), &\text{ if }i=2,
\end{cases}\ 
s_{12}^-(j)(\alpha_i(b))=
\begin{cases}
(t_{b}^{1j2})_*\left(\frac{d}{dz}\right), &\text{ if }i=1,\\
(t_{b}^{2j1})_*\left(\frac{d}{dz}\right), &\text{ if }i=2.
\end{cases}
\end{gather*}
The section $s_{\partial \bar{S}_{12}}$ is defined by 
\[
s_{\partial \bar{S}_{12}}=r_{\bar{S}_{12}}\left(\bigotimes_{j,k}s_{12}(j,k)\otimes\bigotimes_{j=3}^m(s_{12}^+(j)\otimes s_{12}^-(j))\right).
\]

Since any $h\in\mathbf{Z}_d$ fixes $S_{12}$ and $\{S_i\}_{i=3}^m$ pointwise, $q:E\to\bar{E}$ induces an isomorphisms
\[
N(S_{12})^{\otimes d(m-1)(m-2)}\cong N(\bar{S}_{12})^{\otimes(m-1)(m-2)}, \text{ and }N(S_i)^{\otimes d(m-1)(m-2)}\cong N(\bar{S}_i)^{\otimes(m-1)(m-2)}.
\]
By the definition of the sections $s_{\partial S_i}$ and $s_{\partial S_{12}}$, we have 
\begin{align*}
n(s_{\partial S_i},N(S_i)^{\otimes d(m-1)(m-2)})&=n(s_{\partial\bar{S}_i},N(\bar{S}_i)^{\otimes(m-1)(m-2)}),\\
n(s_{\partial S_{12}},N(S_{12})^{\otimes d(m-1)(m-2)})&=n(s_{\partial\bar{S}_{12}},N(\bar{S}_{12})^{\otimes(m-1)(m-2)}).
\end{align*}

First, we will compute the number $n(s_{\partial \bar{S}_i},N(\bar{S}_i)^{\otimes(m-1)(m-2)})$ for $3\le i\le m$. Since $r_{\bar{S}_i}$ induces an isomorphism $T(\bar{E}/\Delta_1)|_{\bar{S}_i}\cong N(\bar{S}_i)$, we have 
\[
n(s_{\partial \bar{S}_i},N(\bar{S}_i)^{\otimes(m-1)(m-2)})=n(\bigotimes_{j,k}s_i(j,k),T(\bar{E}/\Delta_1)|_{\bar{S}_i}^{\otimes(m-1)(m-2)}).
\]
Let $w=x_2/x_1$ be the inhomogeneous coordinate of the second factor of $\bar{E}=\Delta_1\times\mathbf{CP}^1$.
Since the map $t_{b}^{ijk}:\mathbf{CP}^1\to \bar{f}^{-1}(b)$ is written as
\[
t_{b}^{ijk}(z)=\frac{\alpha_k(b)(\alpha_i(b)-\alpha_j(b))z+\alpha_i(b)(\alpha_j(b)-\alpha_k(b))}{(\alpha_i(b)-\alpha_j(b))z+(\alpha_j(b)-\alpha_k(b))},
\]
the vector $(t_{b}^{ijk})_*(d/dz)$ is described as
\begin{equation}\label{section explicit}
(t_{b}^{ijk})_*\left(\frac{d}{dz}\right)=\frac{(\alpha_i(b)-\alpha_j(b))(\alpha_k(b)-\alpha_i(b))}{\alpha_j(b)-\alpha_k(b)}\left(\frac{d}{dw}\right).
\end{equation}

Suppose $j\ge3$ and $k\ge3$. Let $j',k'$ also be integers such that $3\le j'\le m$, $3\le k'\le m$, and $i,j',k'$ are mutually distinct. Let $s_0:\partial\bar{S}_i\to T(\bar{E}/\Delta_i)|_{\partial \bar{S}_i}$ be the zero section. The intersection form on the first homology group of the sphere bundle of $T(\bar{E}/\Delta_1)|_{\partial \bar{S}_i}$ induces that of the $(\mathbf{R}^2-0)$-bundle $T(\bar{E}/\Delta_1)|_{\partial \bar{S}_i}-s_0(\partial \bar{S}_i)$. By the explicit description (\ref{section explicit}), we can calculate the intersection numbers
\begin{gather}
s_i(j,k)\cdot s_i(j',k')=0,\label{calc-intersec1}\\
(s_i(1,k)\otimes s_i(2,k))\cdot s_i(j',k')^{\otimes2}=(s_i(k,1)\otimes s_i(k,2))\cdot s_i(j',k')^{\otimes2}=0,\label{calc-intersec2}\\
(s_i(1,2)\otimes s_i(2,1))\cdot s_i(j',k')^{\otimes2}=1.\label{calc-intersec3}
\end{gather}
The section $s_i(j',k')$ of $T(\bar{E}/\Delta_1)|_{\partial \bar{S}_i}$ can be extended to the nonzero section $\tilde{s}_i(j',k')$ of $T(\bar{E}/\Delta_1)|_{\bar{S}_i}$ defined by $\tilde{s}_i(j',k')(\alpha_i(b))=(t_b^{ij'k'})_*\left(\frac{d}{dz}\right)$ for $3\le i\le m$. Hence a trivialization of the bundle $T(\bar{E}/\Delta_1)|_{\bar{S}_i}$ is given by $s_i(j',k')$. By Lemma \ref{intersection}, we have 
\[
n(s_{\partial \bar{S}_{i}},N(\bar{S}_i)^{\otimes(m-1)(m-2)})=-\left(\bigotimes_{j,k}s_i(j,k)\right)\cdot s_i(j',k')^{\otimes{(m-1)(m-2)}}.
\]
The calculations (\ref{calc-intersec1}), (\ref{calc-intersec2}), and (\ref{calc-intersec3}) show that this is equal to $-1$.

Next, we will compute $n(s_{\partial\bar{S}_{12}},N(\bar{S}_{12})^{\otimes(m-1)(m-2)})$. Suppose $j\ge3$ and $k\ge3$. The section $s_{12}(j,k)$ of $T(\bar{E}/\Delta_1)|_{\partial \bar{S}_{12}}$ can be extended to a nonzero section $\tilde{s}_{12}(j,k)$ of $T(\bar{E}/\Delta_1)|_{\bar{S}_{12}}$ defined by $\tilde{s}_{12}(j,k)(\alpha_i(b))=(t_b^{ijk})_*\left(\frac{d}{dz}\right)$ for $i=1,2$. The section $r_{\bar{S}_{12}}\tilde{s}_{12}(j,k)$ of $N(\bar{S}_{12})$ intersects the zero section $s'_{0}:\bar{S}_{12}\to N(\bar{S}_{12})$ transversely in one point $(0,[0:1])\in \bar{S}_{12}$. Hence we have $n(r_{\bar{S}_{12}}s_{12}(j,k), N(\bar{S}_{12}))=1$. By the explicit description (\ref{section explicit}), We also have
\[
r_{\bar{S}_{12}}s_{12}^\epsilon(j)\cdot r_{\bar{S}_{12}}s_{12}(j,k)=s_{12}^\epsilon(j)\cdot s_{12}(j,k)=-1
\]
on $H_1(N(\bar{S}_{12})-s'_0(\bar{S}_{12}))$ for $\epsilon=+,-$. Lemma \ref{intersection} shows
\begin{equation}\label{section S_i}
n(r_{\bar{S}_{12}}s_{12}^\epsilon(j),N(\bar{S}_{12}))=-r_{\bar{S}_{12}}s_{12}^\epsilon(j)\cdot r_{\bar{S}_{12}}s_{12}(j,k)+n(r_{\bar{S}_{12}}s_{12}(j,k),N(\bar{S}_{12}))=2.
\end{equation}
Thus we have 
\begin{align*}
&n(s_{\partial \bar{S}_{12}}, N(\bar{S}_{12})^{\otimes(m-1)(m-2)})\\
=&\sum_{\begin{subarray}{c}3\le j\le m\\3\le k\le m\\j\ne k\end{subarray}}n(r_{\bar{S}_{12}}s_{12}(j,k),N(\bar{S}_{12}))+\sum_{j=3}^m n(r_{\bar{S}_{12}}s_{12}^{+}(j),N(\bar{S}_{12}))+\sum_{j=3}^m n(r_{\bar{S}_{12}}s_{12}^{-}(j),N(\bar{S}_{12}))\\
=&(m+1)(m-2).
\end{align*}
For any $h\in\mathbf{Z}_d-\{0\}$, the fixed point set $E^h$ is the disjoint union of $S_{12}$ and $\{S_i\}_{i=3}^m$. By the above computations of $n(s_{\partial\bar{S}_i},N(\bar{S}_i)^{\otimes(m-1)(m-2)})$ and $n(s_{\partial \bar{S}_{12}}, N(\bar{S}_{12})^{\otimes(m-1)(m-2)})$, we have 
\begin{align*}
&n(s_{\partial E^h}, N(E^h)^{\otimes d(m-1)(m-2)})\\
=&\sum_{i=3}^{m}n(s_{\partial S_{i}}, N(S_i)^{\otimes d(m-1)(m-2)})+n(s_{\partial S_{12}}, N(S_{12})^{\otimes d(m-1)(m-2)})\\
=&m(m-2).
\end{align*}
Thus the local euler number is 
\[
\chi_{\loc}^{\frac{2h\pi}{d},h}([f,E,\Delta_1])=\frac{m}{d(m-1)}.
\]
Since there are no vertical components, by Theorem \ref{thm:main}, we have
\[
\sigma_{\loc}([f,E,\Delta_1])=-\sum_{h=1}^{d-1}\chi_{\loc}^{\frac{2h\pi}{d},h}([f,E,\Delta_1])\cosec^2\left(\frac{h\pi}{d}\right)+\Sign \bar{E}.
\]
It is known that (for example, see Hirzebruch-Zagier \cite{hirzebruch1974ast} p.178)
\[
\sum_{h=1}^{d-1}\cosec^{2}\left(\frac{h\pi}{d}\right)=\frac{d^{2}-1}{3}.
\]
There is a deformation retraction of $\bar{E}$ onto $\{0\}\times\mathbf{CP}^1$ and its self-intersection number is $0$. Hence we have $\Sign \bar{E}=0$. Thus the local signature and the cobounding function of the pullback of the Meyer cocycle is
\[
\sigma_{\loc}([f,E,\Delta_1])=-\phi(\hat{\sigma}_{ij})=-\frac{(d-1)(d+1)m}{3d(m-1)}.
\]
\end{proof}
\begin{acknowledge}\normalfont
\textrm
The author would like to thank his advisor Nariya Kawazumi for his encouragement and helpful suggestions. He also would like to thank Mikio Furuta for fruitful advices. This research is supported by JSPS Research Fellowships for Young Scientists (20-8182)
\end{acknowledge}

\nocite{namba2003dff} \nocite{matsumoto1996lfg} \nocite{birman1975bla} \nocite{kawazumi1998rhf}
\bibliographystyle{jplain}

\noindent
\textsc{%
Masatoshi Sato\\
Graduate School of Mathematical Sciences, \\
The University of Tokyo \\
3-8-1 Komaba Meguro-ku Tokyo 153-8914, Japan}

\noindent\texttt{%
E-mail address:sato@ms.u-tokyo.ac.jp
}
\end{document}